\documentclass[12pt, a4paper]{amsart}
\usepackage{vmargin}
\usepackage{amsmath, amsthm,thmtools,amssymb}
\usepackage{xypic}
\usepackage[colorlinks=true,citecolor=blue]{hyperref}
\usepackage{epigraph}
\usepackage{graphicx}
\usepackage{mathrsfs}
\usepackage{natbib}
\usepackage{enumitem}
\let\cite=\citep
\usepackage{cleveref}
\newcommand{\be}{\begin{equation}}
\newcommand{\ee}{\end{equation}}
\newcommand{\bes}{\begin{equation*}}
\newcommand{\ees}{\end{equation*}}
\newcommand{\bea}{\begin{eqnarray}}
\newcommand{\eea}{\end{eqnarray}}
\newcommand{\beas}{\begin{eqnarray}}
\newcommand{\eeas}{\end{eqnarray}}
\newcommand{\ben}{\begin{note}}
\newcommand{\een}{\end{note}}
\newcommand{\bexl}{\vskip0.1em\noindent\hrulefill\vskip1em\begin{ExerciseList}}
\newcommand{\eexl}{\end{ExerciseList}\hrulefill}

\newcommand{\bthm}{\begin{theorem}}
\newcommand{\ethm}{\end{theorem}}
\newcommand{\bpro}{\begin{prop}}
\newcommand{\epro}{\end{prop}}
\newcommand{\bcor}{\begin{corollary}}
\newcommand{\ecor}{\end{corollary}}
\newcommand{\bcon}{\begin{conjecture}}
\newcommand{\econ}{\end{conjecture}}
\newcommand{\bp}{\begin{proof}}
\newcommand{\ep}{\end{proof}}
\newcommand{\blem}{\begin{lemma}}
\newcommand{\elem}{\end{lemma}}
\newcommand{\bn}{\begin{note}}
\newcommand{\en}{\end{note}}
\newcommand{\benum}{\begin{enumerate}}
\newcommand{\eenum}{\end{enumerate}}
\newcommand{\bed}{\begin{defn}}
\newcommand{\eed}{\end{defn}}
\newcommand{\brem}{\begin{remark}}
\newcommand{\erem}{\end{remark}}

\newcommand{\btik}{\begin{tikzpicture}\begin{axis}[scale=0.5,axis y line=center, axis x line=middle]}
\newcommand{\etik}{\end{axis}\end{tikzpicture}}

\let\into=\hookrightarrow

\newcommand{\upperRomannumeral}[1]{\uppercase\expandafter{\romannumeral#1}}

\newtheorem{theorem}[equation]{Theorem}      \newtheorem{lemma}[equation]{Lemma}          \newtheorem{corollary}[equation]{Corollary}

\theoremstyle{definition}
\newtheorem{conj}[equation]{Conjecture}

\theoremstyle{definition}
\newtheorem{defn}[equation]{Definition}
\theoremstyle{remark}

\theoremstyle{definition}
\newtheorem{remark}[equation]{Remark}

\numberwithin{equation}{subsection}

\let\into=\hookrightarrow
\let\isom=\simeq

\let\tensor=\otimes

\newcommand{\C}{{\mathbb C}}

\newcommand{\Pic}{{\rm Pic\,}}
\newcommand{\Q}{{\mathbb Q}}

\newcommand{\Spec}{{\rm Spec}}

\newcommand{\Z}{{\mathbb Z}}

\renewcommand{\int}{\operatorname{int}}
\renewcommand{\O}{{\mathcal O}}
\renewcommand{\P}{{\mathbb P}}
\renewcommand{\wp}{{\mathfrak p}}

\newcommand{\benumlab}{\begin{enumerate}[label={\bf(\arabic{*})}]}
\begin{document}

\title{On surfaces satisfying $q=0,p_g=0,c_1^2=9$}\author{Kirti Joshi}\address{Math. department, University of Arizona, 617 N Santa Rita, Tucson
85721-0089, USA.}

\thanks{}\subjclass{}\keywords{}

\newcommand{\fixnumberwithin}[1]{
\numberwithin{equation}{#1}
	\numberwithin{theorem}{#1}
	\numberwithin{conj}{#1}
	\numberwithin{lemma}{#1}
	\numberwithin{proposition}{#1}
	\numberwithin{corollary}{#1}
	\numberwithin{defn}{#1}
	\numberwithin{remark}{#1}
	\numberwithin{rem}{#1}
	\numberwithin{question}{#1}
}

\newcommand{\nws}{\fixnumberwithin{section}}
\newcommand{\nwss}{\fixnumberwithin{subsection}}
\newcommand{\nwsss}{\fixnumberwithin{subsubsection}}
 \fixnumberwithin{section}

\begin{abstract}
I consider the class of surfaces $X$ over algebraically closed fields with numerical invariants given in the title. In characteristic zero, this class contains fake projective planes which were introduced by David Mumford. I prove that in characteristic $p>0$ such surfaces are Hodge-Witt and also ordinary under additional assumptions. In particular, fake projective planes are Hodge-Witt (\Cref{th:hw}). I show that if $X$ is Frobenius split then $X\isom \P^2$ (\Cref{th:f-split}). I also establish a characteristic free characterization of the projective plane using the Nori fundamental group scheme (\Cref{th:simply-connected}). Finally, I show that any fake projective plane over a number field has good ordinary reduction at all but finitely many primes  and in particular fake projective planes exist in positive characteristics (\Cref{th:reduction}).
\end{abstract}
\maketitle

\numberwithin{equation}{section}
\tableofcontents
\section{Introduction}
All surfaces considered in this paper are assumed to be geometrically connected, smooth, projective and minimal surfaces. From \cite{mumford1979}, one has the following  characterization of the projective plane $\P^2$ over complex
numbers: any simply connected smooth, complex, projective
surface with 
\be\label{eq:surf}q=0,p_g=0,c_1^2=9\ee is isomorphic to $\P^2$. 

There is no algebraic proof of this fact. Using Enriques' classification \cite{Mumford1969} and \eqref{eq:surf}, one reduces to ruling out the possibility that the canonical bundle $K$ is ample. As  \cite{mumford1979} notes, this is analytic and goes via \cite{yau77} which shows that any smooth, complex surface with $K$ ample satisfies $c_1^2\leq 3c_3$, with equality holding if and only if $X$ is
uniformized by a cocompact lattice $\Gamma\subset PU(2,1)$ acting on the
unit ball in $\C^2$.  Since \eqref{eq:surf} implies, $c_2^2=3c_2$, it follows from \cite{yau77} that such an $X$ is not simply connected. Thus under the simply connectedness hypothesis, the possibility that $K$ is ample is eliminated and so $K$ is anti-ample and hence of Kodaira dimension $-\infty$ and hence one completes the proof by appealing to Enriques' classification of surfaces.

A surface of general type which satisfies \eqref{eq:surf} is called a
\emph{fake projective plane}.  As was pointed out in \cite{mumford1979}
there are a finite number of fake projective planes and recently
they have all been classified in \cite{prasad2007} (also see \cite{prasad2010}).

This paper studies this question over
fields of characteristic $p>0$. First, I compute several characteristic $p>0$ invariants of a surface satisfying \eqref{eq:surf} and show (in  \Cref{th:hw}) that the slope spectral sequence (\cite[II, 3.1]{Illusie1979a}) of $X$ degenerates at $E_1$ i.e. $X$ is a Hodge-Witt surface. \Cref{th:ordinarity} is substantially more delicate and establishes that if $X$ satisfies Hodge symmetry then $X$ is ordinary in the sense of \cite[Chapitre IV, D\'efinition 4.12]{Illusie1983a}. Next  I prove in \Cref{th:f-split}, that any Frobenius split surface satisfying \eqref{eq:surf} is a projective plane. In  \Cref{th:simply-connected}, I provide a characteristic free version of the characterization of $\P^2$ considered in \cite{mumford1979}. 
In \Cref{th:reduction}, I prove that a fake projective plane defined over a number field has good ordinary reduction at all but finite number of primes. In particular, fake projective planes exist in positive characteristics.

I thank Christian Liedtke and H\'el\`ene Esnault for correspondence.

\section{A class of surfaces}
For a surface $X$, let $q=q(X)$ be its irregularity, $p_g=p_g(X)$ be its geometric genus and let $K$ be its canonical divisor. For any pair of divisors $D_1,D_2$ on $X$, let $D_1\cdot D_2$ be the intersection pairing. Consider the class of surfaces satisfying
\be\label{eq:surf-class}
\begin{aligned}
q(X)&=0\\
p_g(X)&=0,\\ 
K_X^2&=9.
\end{aligned}
\ee

Following \cite{mumford1979}, one defines
\begin{defn}\label{def:fake-proj-plane}
Let $k$ be an algebraically closed field.  A \textit{fake projective plane} over $k$ is a smooth, projective surface $X/k$ such that 
\benumlab
	\item $X$ is of general type, and
	\item  $X$ satisfies \eqref{eq:surf-class}.
\eenum	
\end{defn}

\section{Hodge-Witt and Ordinarity Properties}
Let $k$ be an algebraically closed field of characteristic $p>0$, $W=W(k)$ be the ring of Witt vectors of $k$ and let $H_{cris}^*(X/W)$ denote the crystalline cohomology of $X$ (\cite{Berthelot1976}). Let $H^j(X,W\Omega_X^i)$ denote the de Rham-Witt cohomology of $X$ (\cite{Illusie1979a}). The notion of \textit{Hodge-Witt} and \textit{ordinary variety} are as defined in \cite[Chapitre IV, 4.6, 4.12]{Illusie1983a}. Ordinarity in the sense of \cite[Chapitre IV, 4.12]{Illusie1983a} is stronger than, but implies, the coincidence of Hodge and Newton polygons. In this section I prove the following:
 \bthm\label{th:hw}
Let $k$ be an algebraically closed field of characteristic $p>0$.  Any smooth, projective, minimal surface $X/k$ satisfying \eqref{eq:surf-class} is Hodge-Witt.  In particular, any fake projective plane $X/k$ is Hodge-Witt.
 \ethm
This will be proved later after computing a certain list of invariants of $X$. This is done by means of the following lemmas.
\blem\label{le:chi}
One has $h^{0,1}=0$ and $\chi(\O_X)=1$.
\elem
\bp 
From \cite[2.19]{joshi2020} one has
$$0\leq h^{0,1}-q\leq p_g=0$$
and hence $h^{0,1}=q=0$.
Then as $\chi(\O_X)=1-h^{0,1}+p_g$, one obtains  the second assertion.
\ep

\blem\label{le:c2}
One has $c_2(X)=3$.
\elem
\bp 
From Noether's formula, one has
$$12\cdot\chi(\O_X)=c_1^2+c_2,$$
hence $c_2=3$ follows from $c_1^2=9$ and \Cref{le:chi}.
\ep
My next lemma computes the Hodge-Witt numbers \cite[2.24]{joshi2020} of $X$.
\blem\label{le:hw}
One has 
\benumlab
\item $h_W^{0,1}=h_W^{1,0}=0$.
\item $h_W^{0,2}=h_W^{2,0}=0$.
\item $h_W^{1,1}=1$.
\eenum
\elem
\bp 
For {\bf(1), (2)}, it suffices to prove one of the equalities, the other being implied by Hodge-Witt symmetry \cite[2.26.1]{joshi2020}. By  one has $h_W^{0,1}=\frac{b_1}{2}=q=0$. This proves {\bf(1)}. Again, by  \cite[2.27]{joshi2020}  and \Cref{le:chi}, one has $h_W^{0,2}=\chi(\O_X)-1+\frac{b_1}{2}=0$. This proves {\bf(2)}. Again by \cite[2.27]{joshi2020} and \Cref{le:c2}, one has 
$$h^{1,1}_W=b_1+\frac{5c_2}{6}-\frac{c_1^2}{6}=0+\frac{5\cdot3}{6}-\frac{9}{6}=1.$$
This proves {\bf(3)} and hence all the assertions.
\ep
My next lemma computes the slope numbers \cite[2.22]{joshi2020} of $X$.
\blem\label{le:slope-nums}
One has
\benumlab
\item $m^{0,1}=m^{1,0}=0$.
\item $m^{0,2}=m^{2,0}=0$.
\item $m^{1,1}=1$.
\item $b_2=1$.
\eenum
\elem
\bp
For {\bf(1), (2)}, it suffices to prove one of the equalities, the other being implied by slope-number symmetry \cite[2.25.1]{joshi2020}. One has $m^{0,1}+m^{1,0}=h_W^{0,1}+h_W^{1,0}$ (by \cite[2.25.2]{joshi2020}) and as $m^{i,j}\geq 0$  for all $i,j\geq 0$ by \cite[2.25.1]{joshi2020}, one has $m^{0,1}=m^{1,0}=0$. This proves {\bf(1)}. From \cite[2.25.2]{joshi2020} one has
$$b_2=m^{0,2}+m^{1,1}+m^{2,0}=h_W^{0,2}+h^{1,1}_W+h_W^{2,0}$$
and hence from \Cref{le:hw} and non-negativity of $m^{i,j}$, one obtains {\bf(2), (3), (4)}.
\ep
\bp[Proof of \Cref{th:hw}]
Let $T^{i,j}$ denote the dimension of the domino (defined by \cite[Proposition 2.18]{Illusie1983a}) associated to the differentials $d_1^{i,j}$ in the slope spectral sequence of $X$ (also see \cite[2.23]{joshi2020}). Since $T^{0-1,j+1}=0=T^{0-2,j+2}$ from \cite{Illusie1983a} so from \cite[2.24]{joshi2020} one has
$$h_W^{1,1}=m^{1,1}+T^{0,2}.$$
Thus by \Cref{le:hw} and \Cref{le:slope-nums}, one sees that $T^{0,2}=0$. By \cite[Chapitre I, 2.18.2]{Illusie1983a}, this means  that the de Rham-Witt differential 
$$d:H^2(X,W(\O_X))\to H^2(X,W\Omega_X^1)$$
is zero.  Since $X$ is a surface, by \cite[II, Corollaire 3.14]{Illusie1979a}, all other differential are also zero. This means for $i,j\geq 0$, all the domino numbers $T^{i,j}=0$ i.e. all the differentials in the slope spectral sequence are zero. Hence $X$ is Hodge-Witt by \cite[Chapitre IV, 4.6.2]{Illusie1983a}.
\ep

My next result is substantially more delicate:
\bthm\label{th:ordinarity} 
If Hodge symmetry holds for $X$, then $X$ is ordinary in the sense of \cite[Chapitre IV, D\'efinition 4.12]{Illusie1983a}. Moreover, $H^*_{cris}(X/W)$ is torsion-free.
\ethm
\bp 
Since $X$ is Hodge-Witt, by \cite[Chapitre IV, Th\'eor\`eme 4.5]{Illusie1983a} one has the Hodge-Witt decomposition of finitely generated $W$-modules
\be\label{eq:hw-decomp1}  
H^i_{cris}(X/W)=\bigoplus_{i+j=n} H^i(X,W\Omega_X^j).
\ee

From connectedness of $X$ one sees that 
\be\label{eq:h0} H^0_{cris}(X/W)=W.\ee 
Since Hodge symmetry holds for $X$, one sees that $H^0(X,\Omega_X^1)=H^1(X,\O_X)=0$. This means $H^1_{dR}(X/k)=0$ \cite[Page 64]{oda69}. Hence  from the universal coefficient theorem for crystalline cohomology \cite[II, 4.9.1]{Illusie1979a}, one sees that $$H^1_{cris}(X/W)\tensor_W k=0$$ and
$H^2_{cris}(X/W)$ is torsion free. Hence by Nakayama's lemma, 
\be\label{eq:h1} H^1_{cris}(X/W)=0\ee 
since its reduction modulo $p$ is zero. 

Hence from the Hodge-Witt decomposition \eqref{eq:hw-decomp1}, one has
\be H^1(X,W\Omega^0_X)=H^0(X,W\Omega^1_X)=0.\ee
Since $H^2_{cris}(X/W)$ is torsion-free, and $b_2=1$, one obtains an isomorphism of $W$-modules
\be\label{eq:h2} H^2_{cris}(X/W)\isom W.\ee
Since $X$ is projective, the first Chern class of the ample line bundle is contained in $H^1(X,W\Omega_X^1)$ and hence $H^1(X,W\Omega_X^1)\neq 0$. Hence 
\be 
W\isom H^1(X,W\Omega_X^1).
\ee
This gives (using \eqref{eq:hw-decomp1}) 
\be\label{eq:h20-h02-vanish}  H^2(X,W(\O_X))=0=H^0(X,W\Omega_X^2).\ee
Now I claim that $H^3_{cris}(X/W)=0$. This is quite subtle and does not follow from $H^1_{cris}(X/W)=0$ because Poincar\'e duality pairing is not perfect paring unless one passes to the quotient field of $W$. But at any rate one sees from $H^1_{cris}(X/W)=0$ and Poincare duality that $H^3_{cris}(X/W)$ is torsion. To prove its vanishing one must appeal to the subtler duality theorem of \cite{Ekedahl1984}.

By the Hodge-Witt decomposition \eqref{eq:hw-decomp1} one has
\be H^3_{cris}(X/W)=H^2(X,W(\O_X))\oplus H^1(X,W\Omega^2_X),\ee
and from \eqref{eq:h20-h02-vanish}, one sees that
\be H^3_{cris}(X/W)=H^1(X,W\Omega^2_X).\ee
Hence the latter group is torsion. Now by Ekedahl's Duality Theorem \cite[Theorem 3.5]{Ekedahl1984}, the torsion of $H^1(X,W\Omega^2_X)$ consists of slope-zero,  semi-simple and nilpotent torsion (these terms are defined in \cite[Definition 3.3.13]{Ekedahl1984}). The slope-zero torsion is dual to the slope-zero part of $H^1(X,W(\O_X))$. But this group is zero by \eqref{eq:h20-h02-vanish}. Therefore, $H^1(X,W\Omega_X^2)$ has no slope-zero torsion. The semi-simple torsion of $H^1(X,W\Omega_X^2)$ is dual to the semi-simple torsion of $H^2(X,W(\O_X))$. But again, this group is zero by \eqref{eq:h20-h02-vanish}, and hence $H^1(X,W\Omega_X^2)$ has no semi-simple torsion.  The nilpotent torsion of $H^1(X,W\Omega_X^2)$ is dual to that of $H^3(X,W\Omega_X^{-1})=0$ (all of these assertions are consequences of the subtle nature of Ekedahl's duality theorem cited above). Hence all the three types of the possible torsion of $H^1(X,\Omega_X^2)$ is zero. Hence $H^1(X,W\Omega_X^2)=0$. This proves that 
\be\label{eq:h3} H^3_{cris}(X/W)=0.
\ee 

Now from \eqref{eq:hw-decomp1}, and the fact that $H^i(X,W\Omega_X^i)=0$ if either $i$ or $j$ is an integer not contained in the interval $[0,\dim(X)]$, one sees that
\be 
H^4_{cris}(X/W)=H^2(X,\Omega_X^2).
\ee
Since $b_0=1$, by Poincar\'e duality, one has $b_4=1$, so modulo torsion of $H^2(X,W\Omega_X^2)$, $H^4_{cris}(X/W)\isom W$.  As shown earlier, the torsion in $H^2(X,BW\Omega_X^2)$ can be ruled out by showing that all the three possible types of torsion in $H^2(X,W\Omega_X^2)$ is zero.
This gives 
\be\label{eq:h4} 
H^4_{cris}(X/W)=H^2(X,\Omega_X^2)=W.
\ee
In particular, one sees from \eqref{eq:h4}, \eqref{eq:h3}, \eqref{eq:h2}, \eqref{eq:h1} and \eqref{eq:h0} that $H^*_{cris}(X/W)$ is torsion-free and since the rank of each cohomology group is at most one, for all all $i\geq 0$, $H^i_{cris}(X/W)$  has no fractional slopes of Frobenius.

From \cite[Chapitre IV, Th\'eor\`eme 4.5]{Illusie1983a} one has the decomposition of $W$-modules
\be  
H^j(X,W\Omega^i_X)=H^j(X,ZW\Omega_X^i)\bigoplus H^j(X,BW\Omega_X^{i+1}),
\ee
Hence from the fact $H^*_{cris}(X/W)$ is torsion-free  one sees that $H^j(X,BW\Omega_X^{i+1})$ is $p$-torsion-free for all $0\leq i,j\leq 2=\dim(X)$. On the other hand,  using \cite[Chapitre IV, Th\'eor\`eme 4.5]{Illusie1983a}, for all $i,j\geq 0$ and passing to the quotient field $W[1/p]$ of $W$ one has an isomorphism
\be\label{eq:slope-decomp}  
H^j(X,BW\Omega_X^{i+1})\tensor_WW[1/p] \isom (H^{i+j}_{cris}(X/W)\tensor_WW[1/p])_{]i,i+1[},
\ee
which identifies the cohomology on the left with the portion of $H^{i+j}_{cris}(X/W)\tensor_WW[1/p]$ with slopes of Frobenius lying in the open interval $]i,i+1[$. Since $H^*_{cris}(X/W)$ has no fractional slopes,
 one deduces that the module on the right in \eqref{eq:slope-decomp} is torsion. From the torsion-freeness of $H^*_{cris}(X/W)$,  the module on the right in \eqref{eq:slope-decomp} is free. Thus, one sees that $H^j(X,BW\Omega_X^{i+1})$ is both free and torsion and hence
 one has arrived at
\be  
H^j(X,BW\Omega_X^{i+1})=0 \text{ for all }i\geq 0.
\ee
In other words, by \cite[Chapitre IV, D\'efinition 4.12]{Illusie1983a}, $X$ is ordinary. Along the way, I have established that $H^*_{cris}(X/W)$ is torsion-free. This completes the proof of \Cref{th:ordinarity}.
\ep
\section{The Frobenius split case}
Frobenius splitting property is defined in \cite[Definition 2]{mehta1985}. By \cite[Theorem 2]{mehta1985}, $\P^n$ ($n\geq1$) and hence $\P^2$ is Frobenius split. Hence, the Frobenius splitting hypothesis in the next theorem is quite reasonable.
\bthm\label{th:f-split} 
Let $X/k$ be a smooth, projective, minimal and Frobenius split surface satisfying \eqref{eq:surf-class}. Then $X\isom \P^2$.
\ethm
\bp 
By the hypothesis and  \cite{joshi04b}, $X$ admits a proper, flat lifting to $W_2=W_2(k)$. Hence by \cite{deligne87}, the Hodge-de Rham spectral sequence of $X$ degenerates in dimension $\leq p-1$. As $p_g=0$ so by Serre duality, $H^0(X,\Omega_X^2)=H^2(X,\O_X)=0$, hence for any $p$, the Hodge de Rham sequence degenerates at $E_1$ by \cite{oda69}. 

Let $\kappa(X)$ be the Kodaira dimension of $X$, then Frobenius splitting gives $$H^0(X,K^{1-p})\neq0$$ (\cite[Proposition 6]{mehta1985}) and hence, as observed in the first paragraph of the proof of \cite[Proposition 6.4]{joshi2020}, one has $\kappa(X)\leq 0$. Now I will use the classification given by \cite[Proposition 6.4]{joshi2020}--this needs $p\geq 5$. So assume $p\geq 5$. By the classification of smooth, projective Frobenius split surfaces for $p\geq 5$ given in \cite[Proposition 6.4]{joshi2020}, one has $\kappa(X)\leq 0$ and the following possibilities hold: 
\benumlab
\item $\kappa(X)=-\infty$ and $X$ is one of
\begin{enumerate}
\item rational, or
\item ruled.
\end{enumerate}
\item $\kappa(X)=0$ and $X$ is one of
\begin{enumerate}
\item ordinary K3, or
\item ordinary abelian surface,
\item bielliptic with an ordinary Albanese variety, or
\item ordinary Enriques surface.
\end{enumerate}
\eenum
If {\bf(1)(a)} holds, then $X=\P^2$. So one has to eliminate all the remaining possibilities. If $X$ is ruled over a curve of genus $g$ then by \cite[Chap. V, Corollary 2.11]{Hartshorne1977} one has $K^2=8(1-g) \neq 9$ for any $g\in \Z$. Thus, $X$ is not a ruled surface. Thus {\bf(1)(b)} does not hold.

Now let me consider {\bf(2)}.
If $X$ is a K3 surface then $K^2=0$ and so by \eqref{eq:surf-class}, $X$ cannot be a K3 surface. If $X$ is an abelian or a bielliptic surface then $q\neq 0$ which is impossible by \eqref{eq:surf-class}. So {\bf(2)(b,c)} do not hold. If $X$ is an Enriques' surface then by \cite[II, 7.3.1]{Illusie1979a}, one has $b_2=10$ which contradicts \Cref{le:slope-nums}{\bf(4)}. Thus, the last possibility is eliminated and hence $X\isom \P^2$ as claimed.

Now let me deal with the remaining characteristics $p\leq 3$. If $p\leq 3$, then by \cite{Bombieri1977}, the item {\bf(2)}(c) of the above list has additional members, namely, quasi-elliptic surfaces. By \cite[Proposition, Page 26]{Bombieri1976}, the quasi-elliptic surfaces satisfy $K\equiv0$ and hence $K^2=0$ and the Albanese variety of $X$ is an elliptic curve and hence $q\neq 0$.  Hence, no surface satisfying \eqref{eq:surf-class} is quasi-elliptic. 
Hence, again one is done in case $p\leq 3$. Thus, for all characteristics $p\geq2$, $X\isom \P^2$. Hence, additionally, one has $H^0(X,\Omega_X^1)=0$, Hodge symmetry holds for $X$, and $H^1_{dR}(X/k)=0$ and  $H^*_{cris}(X/W)$ is torsion-free.
\ep

\section{A characteristic free characterization of $\P^2$}
I want to prove the following characteristic free version the above characterization of $\P^2$ obtained in  \cite{mumford1979}.
\bthm\label{th:simply-connected} 
Let $k$ be an algebraically closed field of arbitrary  characteristic (including characteristic zero).
Let $X/k$ be a smooth, projective, minimal surface satisfying \eqref{eq:surf-class}. Assume that 
\benumlab
\item the Nori fundamental group scheme, $\pi_1^{Nori}(X)$, of $X$ is trivial, and 
\item $X$ is not of general type. 
\eenum
Then if $char(k)=0$, then {\bf(1)}$\implies${\bf(2)} and in any characteristic, $X\isom \P^2$.
\ethm
\brem 
Let me remark that by \cite[Chap. II, Proposition 8]{nori1982}, the Nori fundamental group scheme is a birational invariant of smooth, complete varieties.
\erem
\bp 
If $k$ is of characteristic zero, then one may assume $k=\C$. Then from \eqref{eq:surf-class} one sees that $\chi(\O_X)=1$. This together with Noether's formula gives $c_2=3$ (exactly as in \Cref{le:c2}). By \cite[Page 130]{Griffiths1994}, one has the Hodge decomposition:
\be 
H^2(X,\C)=H^{0,2}(X)+H^{1,1}(X)+H^{2,0}(X).
\ee
By \eqref{eq:surf-class}, $p_g=h^{0,2}=\dim H^2(X,\Omega^2_X)=0$ and by Serre duality $0=\dim H^0(X,\Omega^2_X)=\dim H^2(X,\O_X)$. Thus $H^2(X,\C)=H^{1,1}(X)$. Since 
\be c_2=3=\chi_{top}(X)=\sum_{i=0}^4 (-1)^i \dim H^i(X,\C)
\ee
Since $X$ is connected, $\dim H^0(X,\C)=1=H^{4}(X,\C)$ by Poincar\'e duality, and $\dim(H^2(X,\C))\geq 1$ by projectivity of $X$, so one sees that $b_2=\dim H^2(X,\C)=1$ and $H^1(X,\C)=0=H^3(X,\C)$. Using \eqref{eq:surf} and the exponential sequence \cite[Page 51]{Griffiths1994}, one obtains that
\be 
\Pic(X)\isom H^2(X,\Z).
\ee
Since $b_2=1$,  $H^2(X,\Z)$ has rank at most one and all torsion in $H^2(X,\Z)$ lives in $\Pic(X)$ and $H^1(X,\O_X)=0$ says that $\Pic^0(X)=0$. 

Using {\bf(1)} and the surjection of group schemes, $$\pi_1^{Nori}(X)\twoheadrightarrow \pi^{et}_1(X)$$ given by  \cite[Theorem 4.1]{esnault2008} one sees  that \'etale fundamental group scheme $\pi_1^{et}(X)=1$. By \cite[Concluding Remarks (1)]{nori1982} or \cite[Remark 3.2(2)]{esnault2008}, implies that the \'etale fundamental group (again denoted by $\pi_1^{et}(X)$) of $X$ is also trivial i.e $\pi_1^{et}(X)=1$.  Hence the hypothesis {\bf(1)} implies that there are no torsion line bundles on $X$ i.e. $\Pic(X)$ (and hence $H^2(X,\Z)$) is torsion-free. Hence $\Pic(X)=NS(X)=\Z$.

From $NS(X)=\Z$ and \eqref{eq:surf-class} one sees  that $K\neq 0$ and $K$ is either ample (so $X$ is of general type) or anti-ample (so $X$ has Kodaira dimension $\kappa(X)=-\infty$). 

Since $k=\C$, by \cite[Expos\'e XII, Corollaire 5.2]{grothendieck1971a}, one knows that the profinite completion of the topological fundamental group $\pi_1^{top}(X)$ of $X$ is  $\pi_1^{et}(X)$. If $X$ is of general type, then by \cite{mumford1979}, there exists a cocompact lattice $\Gamma \subset PU(2,1)$ such that $X$ is a quotient of a ball by $\Gamma$. Thus, $\Gamma=\pi_1^{top}(X)$ is the topological fundamental group of $X$. Since $\Gamma \subset PU(2,1)$ is a subgroup of a linear group of adjoint type, so by \cite{malcev1940} one sees that $\Gamma$ is residually finite i.e. $\Gamma$ injects into its profinite completion $\hat{\Gamma}=\pi_1^{et}(X)$. Hence, $\pi_1^{Nori}(X)=1$ implies that the \'etale fundamental group is trivial i.e. $\pi_1^{et}(X)=1$ and this implies that the topological fundamental group  $\pi_1^{top}(X)=\Gamma=1$. Thus, $X$ is topologically simply connected,  {\bf(1)}$\implies${\bf(2)} and $\kappa(X)=-\infty$. Hence, by the classification of surfaces with $\kappa(X)=-\infty$ given in \cite{Mumford1969}, one deduces that $X\isom \P^2$.
This completes the proof of $char(k)=0$ case of \Cref{th:simply-connected}.

So from now on, assume that $k$ is algebraically closed of characteristic $p>0$.
I claim that  \Cref{eq:surf-class} and the hypothesis {\bf(1,2)} imply that  the N\'eron-Severi group $NS(X)\isom \Z$. Let $\rho=\rho(X)$ be the rank of $NS(X)$. Since $X$ is projective, $\rho\geq 1$, and from \cite[II, Proposition 5.12]{Illusie1979a} that 
\be 1\leq \rho \leq b_2.
\ee
Since $b_2=1$ by \Cref{le:slope-nums}, one concludes that $\rho=1$ i.e. the rank of the N\'eron-Severi group of $X$ is equal to one. Now one wants to rule out the existence of torsion in $NS(X)$. By \Cref{le:chi}, $H^1(X,\O_X)=0$, so  one sees that $\Pic^0(X)=0$. Since $\pi_1^{Nori}(X)=1$,  one knows that $\Pic^{\tau}(X)=0$ \cite{nori1976} and hence $NS(X)$ is torsion-free. Hence, $NS(X)=\Z$ with an ample line bundle as its generator. In particular,  one sees from this that $K$ is not a torsion line bundle.  Since $K^2=9$,  $K\neq 0$ is either ample or anti-ample. If $K$ is ample, then $X$ is of general type. Hence, $K$ is anti-ample. This means that the Kodaira dimension $\kappa(X)=-\infty$.  Thus, by \cite{Mumford1969}, $X$ is either rational or ruled. If $X$ is ruled over  a smooth, proper curve  $C$ of genus $g\geq1$, then the Albanese variety $Alb(X)=Alb(C)\neq 0$. This contradicts the vanishing of $H^1(X,\O_X)$. Thus, one sees that if $X$ is ruled over $C$, then  $C=\P^1$. In this case $\rho>1$ (\cite[Chap. V, Prop. 2.3]{Hartshorne1977}), which contradicts $\rho=1$. Thus, $X$ is not ruled but has $\kappa(X)=-\infty$ and hence $X\isom \P^2$. This proves the theorem.
\ep

As should be clear from the proof of \Cref{th:simply-connected}, if $X/k$ is a fake projective plane and if $k$ has characteristic zero, then $\pi_1^{Nori}(X)\neq 1$. It is tempting to make the following 
\begin{conj}\label{con:non-triviality}
Suppose $X$ is a fake projective plane over an algebraically closed field of characteristic $p>0$. Then $\pi_1^{Nori}(X)\neq 1$.
\end{conj}

\brem
For some evidence see \Cref{th:reduction}{\bf(2)}.
\erem

\brem 
As was remarked to me by H\'el\`ene Esnault, the proof given above only uses that the assertion that the abelianization of $\pi_1^{Nori}(X)$ is trivial. In characteristic zero, it was shown by \cite[Theorem 10.1]{prasad2007} that if $X$ is a fake projective plane, then the  abelianization $\Pi/[\Pi,\Pi]=\Gamma/[\Gamma,\Gamma]$ of the topological fundamental group is non-trivial i.e. any complex surface satisfying \eqref{eq:surf} with the same integral homology as $\P^2_{\C}$ is isomorphic to $\P^2_\C$. So one could formulate a stronger version of \Cref{con:non-triviality}: if $X$ is any fake projective plane over an algebraically closed field of characteristic $p>0$, then the abelianization  $\pi_1^{Nori}(X)^{ab}\neq 1$.
\erem
\section{Ordinarity of reductions of fake projective planes}
\newcommand{\sX}{\mathscr{X}}
From \cite{mumford1979} or \cite{prasad2007} one knows that fake projective planes are rigid and hence definable over number fields. In particular, this implies, by a standard reduction modulo $p$ argument, the existence of fake projective planes in positive characteristics.  The next theorem asserts ordinarity reductions.
\bthm\label{th:reduction}
	Let $F/\Q$ be a finite extension. Let $X/F$ be a fake projective plane over $F$. Then for all but finitely many non-archimedean primes $\wp$ of $F$, 
	\benumlab
	\item $X$ has good ordinary reduction $X_\wp$ at $\wp$, and 
	\item $\pi_1^{Nori}(X_\wp)\neq 1$, and
	\item $X_\wp$ is a fake projective plane over the algebraic closure of the residue field of $\wp$.
	\eenum
\ethm
\begin{proof}
Let $\O_F\subset F$ be the ring of integers of $F$.  Let $F_\wp$ be the completion of $F$ at some prime ideal $\wp\subset \O_F$. Let $\O_{F_\wp}$ be the ring of integers of $F_\wp$. Let $\O_{F_\wp}^{sh}$ be the strict Henselization of $\O_{F_\wp}$. Choose a proper, flat scheme $\sX/\O_F$ of such that $\sX\times_{\Spec(\O_F)} \Spec(F)\isom X$.  Then there exists a closed subset $V(I)\subset \Spec(\O_F)$ such that for all $\wp\not\in V(I)$,  one has
\benumlab
\item $\sX_\wp=\sX\times_{\Spec(\O_F)} \Spec(\O_{F_\wp}^{sh})$ is proper and smooth over $\O_{F_\wp}^{sh}$, and
\item the de Rham cohomology $H^i_{dR}(\sX_\wp/\O_{F_\wp}^{sh})$ is a free $\O_{F_\wp}^{sh}$-module, and
\item the special fiber $X_\wp$ of $\sX_\wp$ satisfies \eqref{eq:surf}, and
\item if $k(\wp)$ is the residue field $\O_{F_\wp}^{sh}/\wp$ of $\wp$, then $H^i_{dR}(\sX_\wp/\O_{F_\wp}^{sh})\tensor_{\O_{F_\wp}}k(\wp)\isom H^i_{dR}(X_\wp/k(\wp))$.
\eenum
One works with the strict Henselization $\O_{F_\wp}^{sh}$ instead of $\O_{F_\wp}$ because the residue field of $\O_{F_\wp}^{sh}$ is algebraically closed. This is needed because of the requirement in the preceding results that the field $k$ is algebraically closed. Then by {\bf(2,3,4)}, $\sX_\wp$ also satisfies Hodge symmetry and hence by \Cref{th:ordinarity}, $X_\wp$ is ordinary. This proves {\bf(1)}.

Choose an embedding $F\into \C$. Since $\pi_1^{top}(X/\C)\neq 1$ by \cite{mumford1979}, one sees (as shown earlier) that $\pi_1^{et}(X/\C)\neq 1$. Hence by \cite{nori1976} or \cite{esnault2008}, $\pi_1^{Nori}(X/\C)\neq 1$. By  \cite[Expos\'e X, Corollaire 3.9]{grothendieck1971a}, one sees that the prime-to-$p$ part $\pi_1^{et}(X_\wp/k(\wp))^{(p)}$ is isomorphic to the prime-to-$p$ part $\pi_1^{et}(X/\C)^{(p)}$ of $\pi_1^{et}(X/\C)$ and as the latter is non-trivial, one sees that $\pi_1^{et}(X_\wp)\neq 1$. Using this and \cite[Theorem 4.1 and Remark 3.2(2)]{esnault2008}, one sees that $\pi_1^{Nori}(X_\wp/k(\wp))\neq 1$. This completes the proof of the assertion {\bf(2)}. Finally note that $X_\wp$ satisfies \eqref{eq:surf} and is of general type, hence by \Cref{def:fake-proj-plane}, $X_\wp$ is a fake projective plane over $k(\wp)$. This proves {\bf(3)}.
\end{proof}

\brem
\Cref{th:reduction} establishes the existence of fake projective planes in all but finitely many characteristics. Using \cite{prasad2007,prasad2010} and allowing $F$ to vary, one can expect to establish the existence of fake projective planes in all characteristics.
\erem

\bibliographystyle{plainnat}

\end{document}